\documentclass[11pt]{article}

\usepackage{import} 
\usepackage[left=2.5cm,right=2.5cm,top=2.5cm,bottom=2.5cm]{geometry}

\usepackage{subfiles}
\usepackage{appendix}

\usepackage[natbib, backend=bibtex, 
            style=authoryear,
            url=false, 
            ]{biblatex}

\usepackage{microtype} 

\usepackage[hyperfirst=true]{glossaries}

\usepackage{amssymb, amsthm, mathtools}
\usepackage{dsfont} 
\usepackage[ruled,linesnumbered]{algorithm2e}

\usepackage{booktabs} 
\usepackage{multirow} 
\usepackage{diagbox} 

\usepackage{graphicx} 
\usepackage[dvipsnames]{xcolor} 

\usepackage{caption} 
\usepackage{subcaption}
\usepackage{wrapfig} 

\usepackage{silence} 

\usepackage[breaklinks, colorlinks]{hyperref} 
\hypersetup{colorlinks, citecolor=blue, urlcolor=black}
\usepackage[noabbrev]{cleveref} 

\usepackage{minitoc}

\newcommand{\myBullet}{$\bullet~$}

\AtBeginEnvironment{subappendices}{%
    \section*{Appendix}
    \addcontentsline{toc}{section}{Appendices}
}

\newcommand{\mbb}{\mathbb} 
\newcommand{\mc}{\mathcal} 
\newcommand{\mbf}{\mathbf} 

\newcommand{\un}{\mathds{1}} 

\newcommand{\under}[2]{\underset{#2}{#1}}
\newcommand{\som}[3]{\sum_{#1=#2} ^{#3} }
\newcommand{\pro}[3]{\prod_{#1=#2} ^{#3} }

\newcommand{\pen}{\mathrm{pen}}
\newcommand{\regularization}{\lambda}

\DeclareMathOperator{\prox}{Prox}
\DeclareMathOperator*{\argmin}{arg\,min}
\DeclareMathOperator*{\argmax}{arg\,max}

\newcommand{\setr}{\mathbb R}

\newcommand{\sets}{\mathbb S}
\newcommand{\setsppr}[1]{\sets^{++}_{#1}}

\newcommand{\abs}[1]{\left\vert #1 \right\vert}
\newcommand{\norm}[1]{\left\|#1\right\|}

\newcommand{\stackeq}[2][rrrrrrrrr]{\left\{\begin{array}{#1} #2 \end{array}\right.}

\newcommand{\memfct}{m}
\newcommand{\marg}{g}

\newcommand{\hazard}{h}

\newcommand{\para}{\theta} 
\newcommand{\paraspace}{\Theta} 

\newcommand{\nobs}{N} 
\newcommand{\nrep}{J} 

\newcommand{\event}{T} 
\newcommand{\obs}{Y} 

\newcommand{\lat}{Z} 

\newcommand{\varlat}{\Omega} 
\newcommand{\meanlat}{\mu} 

\newcommand{\covrep}{t} 
\newcommand{\eps}{\varepsilon} 
\newcommand{\vareps}{\Sigma} 

\newcommand{\covgd}{X} 
\newcommand{\dimcovgd}{p} 
\newcommand{\ldp}{\tau} 

\newcommand{\iter}{k}

\newcommand{\step}{\gamma}

\newcommand{\Prec}{P}

\theoremstyle{plain}

\theoremstyle{definition}

\theoremstyle{remark}

\newtheoremstyle{exampstyle}
    {\topsep}
    {\topsep}
    {\itshape}
    {}
    {\bfseries}
    {.}
    {.5em}
    {}

\theoremstyle{exampstyle}

\theoremstyle{remark}

\SetKwInOut{Require}{Require}
\SetKw{Initialize}{Initialize}
\SetKw{Return}{Return}

\SetKwFor{Tq}{tant que}{faire}{fin}
\SetKwInput{Entree}{Entrées}
\SetKw{Retour}{retourner}

\newcommand{\codename}[1]{\textsc{#1}\xspace}

\newcommand{\NLJMLASSO}{\codename{NLJM-LASSO}}

\newenvironment{figureht}{
    \begin{figure}[!ht]\centering
	}{\end{figure}}
		
\newenvironment{tableht}{
    \begin{table}[!ht] \centering
	}{\end{table}}

\usepackage[left=2.5cm,right=2.5cm,top=2.5cm,bottom=2.5cm]{geometry}

\usepackage[figuresright]{rotating}

\newcommand{\hatLASSO}[1][\theta]{\hat{\mbf{#1}}_{\mbf{pen}}}
\newcommand{\hatMLE}[1][\theta]{\hat{\mbf{\theta}}_{\mbf{MLE}}}

\def\geno{i}
\def\nobs{N}

\def\nrep{J}

\def\hazard{h}
\def\regularization{\lambda}
\def\memfct{m}

\def\varobsValue{{\boldsymbol{\mc D}}}

\def\varcov{X}
\newcommand{\Lmarg}[1][]{\mc L_{marg#1}}

\newcommand{\varobs}[1][]{\varobsValue#1}
\renewcommand{\varlat}[1][]{\varlatValue#1}

\title{Estimation and variable selection in high dimension in a causal joint model of survival times and longitudinal outcomes with random effects.}

\newcommand{\emailx}[1]{ \it{email:}\href{mailto:#1}{#1}}

\author
{
    Antoine Caillebotte\emailx{caillebotte.antoine@inrae.fr} \\
    Université Paris-Saclay, INRAE, UR MaIAGE, UMR GQE-Moulon, France
    \and
    Estelle Kuhn\emailx{estelle.kuhn@inrae.fr} \\
    Université Paris-Saclay, INRAE, UR MaIAGE, France,
    \and
    Sarah Lemler\emailx{sarah.lemler@centralesupelec.fr} \\
    Université Paris-Saclay, Laboratoire MICS, France,
}

\newenvironment{keywords}{Key words:\qquad}{}
\renewenvironment{abstract}{Abstract:\qquad}{}

\begin{document}




\label{firstpage}

    \maketitle

    \begin{abstract} 
        We consider a joint survival and mixed-effects model to explain the survival time from longitudinal data and high-dimensional covariates in a population. The longitudinal data is modeled using a non linear mixed-effects model to account for the inter-individual variability  in the population. The corresponding  regression function serves as a link function incorporated into the survival  model. In that way, the longitudinal data is related to the survival time.  We consider a Cox model that  takes into account  both  high-dimensional covariates and the link function. There are two  main objectives: first,  identify the relevant covariates that contribute to explaining survival time, and second,  estimate all unknown parameters of the joint model. For the first objective,  we consider the estimate defined by maximizing the marginal log-likelihood regularized with a $\ell_1$-penalty term.
        To tackle the optimization problem, we implement an adaptive stochastic gradient to handle the latent variables of the non linear mixed-effects model associated with a proximal operator to manage the non-differentiability of the penalty.  We rely on an eBIC model choice criterion to select an optimal value for the regularization parameter.  Once the relevant covariates are selected, we re-estimate the parameters in the reduced model by maximizing the likelihood using an adaptive stochastic gradient descent.
              We provide relevant simulations that showcase the performance of the proposed variable selection and  parameter estimation method in the joint model. We investigate the effect of censoring and of the presence of correlation between the individual parameters in the mixed model.
    \end{abstract}

    \begin{keywords}
        Joint model, nonlinear mixed-effects model, Cox model, high dimension covariates, penalized estimation, adaptive stochastic gradient, proximal operator
    \end{keywords}

    \newpage

\section{Introduction}

    A very current issue in many fields is the need for a better understanding of the interactions between dependent dynamic phenomena. For example, in medicine, this may involve the dynamics of a patient's tumors in oncology and the effects of anti-cancer treatments administered to the patient. Another example in plant science is the dynamics of plant development in a plot and the spread of an epidemic disease or pests in that plot. The considered phenomena are often complex, both in terms of their modes of interaction and their temporal and spatial dynamics. Moreover, these phenomena are frequently observed in populations of heterogeneous or structured individuals, such as patients or plants.

    Mathematical modeling has proven to be a powerful tool for understanding the interactions between multiple dynamic phenomena. It also allows for considering variability present in the observed population of individuals. Joint modeling of several phenomena has demonstrated its effectiveness in several fields, including medicine, pharmacology, and biology (\cite{kerioui_modelling_2022}). A particular case of joint models concerns the simultaneous modeling of longitudinal data and survival data observed on the same individual based on a causal relation. In this type of joint model, longitudinal data are often modeled by a mixed-effects model (\cite{pinheiro_mixed-effects_2000, davidian_nonlinear_1995}), and survival data by a survival model such as the Cox model (\cite{cox_regression_1972}). The latter allows for modeling the instantaneous risk of the survival variable as a function of the covariates. It is also possible to include longitudinal data modeling
in a causal way     as a covariate in the Cox model via a link function. Several authors have studied this model (\cite{Tsiatis97}, \cite{rizopoulos_joint_2012}, \cite{krol2016joint}). 
    
   Inference by maximum likelihood is not direct due to the presence of latent variables in the mixed-effects model and can be performed using numerical tools as for example  Expectation Maximization (EM) like algorithms (\cite{Tsiatis97}, \cite{guedj2011joint}, \cite{rizopoulos2010jm}, \cite{desmee2017using}). The EM-type algorithms, such as the classical Stochastic Approximation Expectation Maximization (SAEM), are the most classical approaches for inferring parameters in the presence of latent variables. They have been developed for estimation in general latent variable models. They are straightforward to implement in the context of a curved exponential family based on sufficient statistics of the model. Moreover, theoretical convergence results have been established in this context.
    However, when the model does not belong to the exponential family, which is the case for such causal joint models,  the methodology is not generic in practice, and the theoretical results fail.    Some exponentialization tricks have then been proposed to face this restrictive assumption of the curved exponential family. For instance, one of them consists of considering some unknown parameters as random population variables into a modified model.  However, \cite{debavelaere2021} have shown that, in general, the parameter returned by the SAEM algorithm on the modified model is not a maximum likelihood of the initial model, and they have suggested the use of this exponentialization trick with variances of the new random population variables  decreasing as the iterations of the algorithm progress. This approach also has limitations in practice due to complex algorithmic settings and tuning.  Besides gradient-based methods are another type of approach, often omitted for estimating parameters in latent models \cite{cappe2005springer}. Recently, \cite{Calimero} suggested using a preconditioned stochastic gradient algorithm to deal with parameter estimation in the presence of latent variables in very general latent variables models. This approach is particularly interesting when considering  models that do not belong to the exponential family, as for causal joint models since theoretical convergence properties are established under general assumptions. Moreover the authors  showed that this algorithm performs well for the nonlinear logistic growth mixed-effects model, which can be used to represent some longitudinal data. Note that Bayesian numerical methods have also been proposed in parallel to estimate the parameters of a joint model. (\cite{r2003bayesian}, \cite{JSSv072i07}, \cite{kerioui_modelling_2022}).

    Besides, in many applications,  today's technological means allow for collecting high-dimensional explanatory covariates. These may include, for example, genetic markers or omics data. In addition to the wealth of information provided by these covariates, they also generate difficulties in the statistical analysis of models as it is necessary to adapt statistical and numerical approaches to their high dimensionality. One possible approach to face this high dimensional setting is to consider a parsimonious assumption, for example by assuming that only a small subset of covariates are relevant to explain the variability observed on the survival times. In practice one can consider a regularized estimator, such as the LASSO one (\cite{he_simultaneous_2015}, \cite{yi2022simultaneous}), and adapted numerical methods, such as  proximal  algorithms (\cite{achab_learning_2017}, \cite{fort_stochastic_2017}).

    In this paper, we consider a causal joint model that combines through a link function  a Cox model for the survival times and a nonlinear mixed-effects model for longitudinal data, including covariates of high dimension. The individual variability in the population is taking into account both in the survival part through the individual covariates and in the longitudinal part through the individual parameter of the mixed model. Our work aims to capitialize on the population approach to  efficiently select the relevant covariates among the high-dimensional covariates involved in the Cox model part of the joint model based on the whole dataset and then to estimate in the reduced model the unknown parameters. For that purpose, we propose a procedure combining a variable selection step based on a $\ell_1$-penalized maximum likelihood estimate for the regression parameter of the Cox model with a standard  maximum likelihood estimation step in the reduced model. To calculate the penalized estimate in practice, we develop an algorithm combining an adaptive stochastic gradient descent to deal with the presence of individual parameters in the mixed model part which are latent variables and a weighted proximal algorithm to handle the non-differentiability of the $\ell_1$ penalty term used for variable selection in the Cox model part. The proposed algorithm is easy to implement in general joint models, in particular it is  not require   that the model density belongs to the curved exponential family.

    The paper is organized as follows. In Section \ref{chap4:sec:model}, we present the causal joint model constructed from  a Cox model for survival data, with high-dimensional covariates and a link function related to a non linear mixed-effects model for longitudinal data; we detail some examples. In Section \ref{sec:algo}, we present the proposed inference method based on: i) a variable selection procedure composed of the evaluation of the $\ell_1$-penalized maximum likelihood estimate for several value of the regularization parameter on a grid and of the choice of this regularization parameter using the eBIC model choice criterion; ii) the re-estimation of the parameters in the reduced model once the relevant covariates have been selected. The fourth section details    the numerical steps, in particular the adaptive proximal weighted   stochastic gradient algorithm used to target the $\ell_1$-penalized estimate on the grid. Finally, we illustrate through a simulation study the performance of the proposed methodology in Section \ref{sec:simu}. The paper ends with a conclusion.


\section{Joint causal model for survival times and longitudinal outcomes   }\label{chap4:sec:model}

 Let us begin by introducing some notation.    We consider a population of $\nobs$ individuals. For each individual $\geno$, we denote by $\mbf T_i^*$  the survival time corresponding to the duration until the occurrence of an event of interest.  We also observe for each individual $\geno$ longitudinal data, which are repeated observations over times. We denote by $\mbf Y_{i,j}$ the $j$th observation of the  $i$th individual for  $i \in \{ 1, \dots,  \nobs \}$ and $j \in  \{ 1, \dots,  \nrep \}$. Note that we consider in our work the same number of repeated observations per individual; however it can be easily  generalized to the case where there are different numbers of longitudinal observations for each individual in the population. The following sections describe the considered joint model.
    
    \subsection{Survival Model}

 For  $i \in \{ 1, \dots,  \nobs \}$   the survival time $\mbf T_i^*$ of individual $i$ represents the time between a fixed initial moment and the occurrence of an event of interest and is modeled through a random variable. To characterize the distribution of $\mbf T_i^*$, we consider the hazard function defined by:
    
        \begin{equation}\label{eq:hazard}
            \hazard_i(t) := \lim_{\Delta t\to 0} \frac{\mbb P (t \leq \mbf T_i^* < t + \Delta t | \mbf T_i^* \geq t) }{\Delta t} 
            ; \forall t \geq 0.
        \end{equation} 

    The Cox model (see \cite{cox_regression_1972}) is one of the most classical models in survival analysis. It allows us to relate the hazard function of the survival time $\mbf T_i^*$ to  covariates $\varcov_i\in\setr^p$, with $p$ being the number of covariates, characterizing individual $i$.  In our approach, we will consider the high-dimensional setting with many covariates so that $p$ is very large with respect to $N$. The hazard of the Cox model for individual $i$ is written as follows:
    
        \begin{equation}\label{eq:cox}
            h(t|\varcov_i) = h_0(t) \exp(\beta^T \varcov_i)
            ; \forall t \geq 0,
        \end{equation}
        
    \noindent with $\beta \in\setr^p$ a regression parameter and $\hazard_0$ the baseline hazard function that characterizes a common behavior in the observed population. In the sequel, we  consider a parametric baseline function denoted by  $\hazard_{\theta_b }$ where $\theta_b$ are its parameters. Therefore, the Cox model's unknown parameters are $\beta$ and $\theta_b$.  
    
  In many settings  the survival time is not observed directly due to censorship (see \cite{cox_regression_1972}).  Thus let's suppose that $\mbf T_i^*$ is censored on the right, and let's denote by $C_i$ the censoring time. We then observe $\mbf T_i= \min(\mbf T_i^*, C_i)$ and $\delta_i = \un_{\mbf T_i^* \leq C_i}$ respectively the censored survival time and the censoring indicator.

    An intrinsic assumption inherent to the Cox model is that the  hazards are proportional, meaning that the ratio between the hazard of two individuals only depends on a function of the covariates of these two individuals and does not depend on the time. This assumption is quite restrictive, in particular many dataset do not satisfy it in practive. For this reason, several extensions have been developed to relax this assumption. It is possible for example to include time-dependent covariates or to consider a random effect as in the frailty models. When additional longitudinal data related  to the survival data are available, another way to release this assumption is to  explain some of the survival time variability using the longitudinal data dynamic, by  modeling them with a nonlinear mixed-effects model and including a causal effect through a link into the Cox  model. 
    
    Let us present the mixed-effects model before explaining the integration of this new component into the Cox model.
        
    \subsection{Non linear mixed-effects model}

    The longitudinal data are observed $\nrep$ times for each individual $i \in \{ 1, \dots,  \nobs \}$. Let us  denoted by $\mbf Y_{i,j}$ the $j$th observation of the $i$th individual for $j \in \{ 1, \dots,  \nrep \}$ and $i \in \{ 1, \dots,  \nobs \}$. We model this longitudinal observation using a non linear function $\memfct$ that depends on individual parameters represented by the latent variable $\lat_i$  as follows:
    
    \begin{equation}\label{eq:nlmem}
        \stackeq[clr]{
            \mbf Y_{i,j} &= \memfct(t_j; \lat_i) + \varepsilon_{i,j},  & \forall 1 \leq i \leq \nobs, 1 \leq j \leq \nrep
            \\ 
            \lat_i &\under\sim{i.i.d.} \mc N(\mu, \Gamma)  &\forall 1 \leq i \leq \nobs \\
           \varepsilon_{i,j} &\under\sim{i.i.d.} \mc N(0,\sigma^2) &\forall 1 \leq i \leq \nobs, 1 \leq j \leq \nrep
        }
    \end{equation}
    \noindent where, $t_j$ is the $j$th observation time, and $\varepsilon_{i,j}$ is an additive noise assumed centered Gaussian with unknown variance $\sigma^2$. The latent variable $\varlat_i$ describes the inter-individual variability of the population. It is assumed that $\varlat_i$ follows a Gaussian distribution with unknown expectation  $\mu$ and covariance matrix  $\Gamma$. 
    We assume that for each $i$ the variables $\varlat_i$ and $(\varepsilon_{i,j})_j$ are independent. We assume also that the variables $(\varlat_i)_i$ are independent, as well as the variables $(\varepsilon_{i,j})_{i,j}$.
    The unknown parameters of the non linear mixed effects model are therefore $\mu, \Gamma$  and $ \sigma^2$.

    Let us introduce in the following the link function, which will combine the two previous models by modeling the  influence of the dynamic of the longitudinal observation of the hazard function.

    \subsection{Joint  causal model}

    We introduce the link function denoted by $f$, which is linearly parameterized with $\alpha$ and depends on the history of the true unobserved longitudinal dynamic. Let us define the set   $\mc M(t;\lat_i) = \{ \memfct(s;\lat_i) |\forall s, 0 \leq s < t \}$ that contains the values of the longitudinal dynamic up to time $t$. Several forms for $f$ can be taken, an example of which is given below. We assume that the hazard of the survival time of individual $i$ is related to the longitudinal  dynamic through the link function $f$ as follows, $\forall 1 \leq i \leq \nobs, \forall t \geq 0$ :
    
        \begin{equation}\label{eq:joint}
                \hazard(t|\mc M(t,\lat_i), \varcov_i) 
            = \hazard_{\theta_b }(t) \exp \bigr(\beta^T \varcov_i+  f(\alpha,\mc M(t;\lat_i)) \bigl).
        \end{equation}

    The parameter $\alpha$ represents the influence of the longitudinal dynamic on the survival data.   The unknown parameters for the joint model include the parameters of the Cox model and those of the non linear mixed-effects model, as well as the link function parameter of the joint model. We denote by $\theta = (\theta_b, \beta, \alpha, \mu, \Gamma,\sigma^2)$ the vector of unknown parameters taking value in  the parameter space denoted by  $\Theta\subset \setr^d$.
    
\section{Variable selection and estimation procedure}\label{sec:algo}
    In this section, we propose a method for  selecting the relevant covariates and  estimating the model parameters.

    \subsection{Inference in latent variables models}
    We  consider the maximum likelihood estimator to infer the joint model parameters.
    Since the random effects of the mixed model part are non observed, it is not possible to evaluate the full complete likelihood of the joint model depending on both observations and latent variables. Therefore as done usually in 
          the context of latent variable models, we consider  the marginal likelihood, denoted by $\Lmarg$   obtained by integrating the full complete likelihood over the latent variables. Let $\varobs_i= (\mbf Y_i, \mbf T_i, \delta_i)$ be the  observed variables for $1 \leq i\leq N$ and $\varobs=(\varobs_i)_{1 \leq i\leq N}$. The marginal likelihood writes:

        \begin{eqnarray}
            \Lmarg(\theta;\varobs)  &=& \pro i1n\int p_\theta(\varobs[_i],\lat_i) d\varlat  [_i]  \notag
           =\pro i1n\int  p_\theta(\varobs[_i]|\lat_i)p_\theta(\lat_i) d\lat_i \label{eq:Lcomp},
        \end{eqnarray}
        
        \noindent where $p_\theta(\varobs,\varlat), p_\theta(\varobs|\varlat), p_\theta(\varlat)$ are respectively  the density of the pair $(\varobs,\varlat)$, the density of $\varobs$ conditionally to $\varlat$, and the density of $\varlat$. We assume that the variables $(\mbf Y_i)$ and $(\mbf T_i)$ are  independent conditionally to the variables $(\varlat_i)$. We assume also that the censoring times $(C_i)$ and $(\mbf T_i)$ are  independent conditionally to the variables $(\varlat_i)$.  Thus we get the following expression for the likelihood:

        \begin{align*}
            \Lmarg (\theta;\varobs) =& \pro i1N\int   p_\theta(\mbf Y_i|\lat_i)
       p_\theta(\mbf T_i, \delta_i|\lat_i) p_\theta(\lat_i) d\lat_i
            \\ =& \pro i1N\int \Bigg(\left\{ \pro j1J \Phi(Y_{i,j}; m(t_j,\lat_i),\sigma^2)\right\}
          \times  \left(h(T_i|\mc M(T_j,\lat_i), \varcov_i)\right)^{\delta_i} 
            \\ &\times \exp\left(-\int_0^{T_i} h(s|\mc M(s,\lat_i), \varcov_i)ds \right)
         \times \Phi(\lat_i; \mu, \Gamma)\Bigg) d\lat_i,
        \end{align*}

        \noindent where $\Phi(.; m,v)$ is the density of the normal distribution with expectation  $m$ and covariance $v$.

Model parameters are typically estimated  using the maximum likelihood estimator, defined as:
\begin{equation}\label{eq:MLE_estimator}
            \hatMLE(\varobs)= \argmax_{\theta\in \Theta} \log \Lmarg(\theta;\varobs) ,
        \end{equation}

        \noindent where $\Theta$ denotes the parameter space.

However, in practice, it is difficult to directly compute this estimate, which does often not have an analytical form in  latent variable models. Therefore, one has to use numerical methods. When dealing with latent variables, traditional approaches for computing this estimate  often rely on Expectation-Maximization (EM) algorithms (see \cite{ng2012algorithm}). A key limitation of these methods is that they are primarily suited to models within the curved exponential family, both from theoretical and practical points of view.  Indeed theoretical results are established ony for models of the curved exponential family. Moreover implementation outside the exponential family is complex and require to store numerous quantities in memory. To overcome  this limitation, a common strategy involves applying the exponentiation trick, which consists of modeling some parameters as Gaussian random variables to construct an augmented model belonging to the exponential family. Nonetheless, this approach also has limitations, as discussed by \cite{debavelaere2021}. Another way to address this issue is to implement  stochastic gradient descent algorithms to solve the optimization problem using  the Fisher identity (see \cite{cappe2005springer}).  In particular a preconditioned stochastic gradient descent algorithm for parameter estimation in latent variable models was  recently proposed by \cite{Calimero}. This method is more generally applicable and theoretical convergence results are stated without requiring that  the model belongs to the curved exponential family.

    \subsection{Definition of the penalized estimator for variable selection }

        In the presence of covariates of high dimension, classical estimation strategies are no more efficient, so we introduce a penalty term to regularize the estimation and perform the variable selection.  We borrow the idea of the LASSO (Least Absolute Shrinkage and Selection Operator) procedure,  developed initially for linear regression models (see \cite{tibshirani_regression_1996}) and extended in several models, including the Cox model (see \cite{tibshirani_lasso_1997}). This method enables  to handle high-dimensional covariates and to select a subset of explanatory covariates from a large collection.
        Assuming only a few variables are relevant to be included in the model among the whole collection,  we consider a $\ell_1$ penalized maximum likelihood estimator to select them.  We emphasize that we consider a  penalty term, which only depends on the parameter $\beta$:

        $$\pen_\regularization(\theta) = \regularization \norm{\beta}_1 = \regularization \som k1p \abs{\beta_k},$$

        \noindent where $\regularization$ is a positive real called regularization parameter.  Let us define for each value of  $\regularization$ in a grid $\Lambda$ the penalized maximum likelihood estimator by:

        \begin{equation}\label{eq:penalized_estimator}
            \hatLASSO^\regularization (\varobs)= \argmax_{\theta\in \Theta} \left(\log \Lmarg(\theta;\varobs) - \pen_\regularization(\theta)\right),
        \end{equation}

        \noindent where $\Theta$ denotes the parameter space and $\regularization$ is a positive parameter. The larger the value of $\regularization$, the more $\beta$ will be constrained to have zero components. Conversely, the smaller the value of $\regularization$, the freer the components of $\beta$ will be. It is customary to determine the value of $\regularization$ using cross-validation (\cite{tibshirani_regression_1996}), in particular for prediction objective. But here, without any predictive context, we choose to  fix an optimal value for the  regularization parameter based on a model criterion as detailed in the next section.

        \subsection{Regularization path procedure}
    
        To select only the most explanatory covariates, it's important to choose a well-balanced value for the regularization parameter $\regularization$. We choose an optimal parameter by minimizing the extended Bayesian Information Criterion (eBIC) (see \cite{eBIC_chen_chen_2008}).  Following the study of  \cite{delattre_mem_BIC_2014}, we penalize the number of degrees of freedom by the logarithm  of the total number of observations $\nobs\times\nrep$. 
        We consider  a grid  $\Lambda$ of values for the regularization parameter $\regularization$. We conduct then the following inference methodology \label{algo:metho}: 
            
        \begin{enumerate}
            \item[i)] For all $\regularization\in\Lambda$ repeat the following steps:
            \begin{itemize}
\item            calculate  $\hatLASSO^\regularization (\varobs)= \argmax_{\theta\in \Theta} \left(\log \Lmarg(\theta;\varobs) - \pen_\regularization(\theta)\right),
       $

                \item  deduce the associated support $\hat S_\regularization=\{j\in\{1,\dots,p\}, \hat\beta_{\regularization,j}^{PEN}\neq 0\}$
                \item calculate  $
                \hat\theta^{\text{MLE}}_\regularization(\varobs)= \argmax_{\theta\in \paraspace} \left\{\log \marg_\regularization(\theta;\varobs) \right\},
                $ where  $\marg_\regularization$ is the marginal likelihood in the  model reduced to the support  $ \hat S_\regularization$.
                \item calculate the eBIC criterion:
                \begin{equation}\label{eq:eBIC}
                    \text{eBIC}(\regularization)=-2\log \marg_\regularization(\hat\theta^{MLE}_\regularization; \obs)+|\hat S_\regularization|\log(\nobs\nrep) +  2\log\left(\binom{\dimcovgd}{|\hat S_\regularization|}\right)
                \end{equation}
                \end{itemize}

 \item[ii)] Select the parameter $\regularization$ that minimizes the eBIC:
                \[\hat\regularization=\argmin_{\regularization\in\Lambda}\text{eBIC}(\regularization),
                \]
                and consider the final estimator  defined by $\hat\theta^{MLE}_{\hat\regularization}$.
        
        \end{enumerate}

   \section{Numerical methods for inference procedure}

To implement our variable selection and re-estimation procedures in practice, we must be able to compute both the maximum likelihood estimate in a reduced model and the penalized estimate across a range of regularization parameter values, with both efficiency and accuracy. In both cases, we face the challenge of integrating over latent variables, which typically precludes closed-form expressions for the optimization criteria. As a result, we turn to numerical methods to solve these maximization problems. Specifically, we propose using an adaptive stochastic gradient algorithm for the first objective and a proximal extended variant for the second.

        \subsection{Adaptive Stochastic Gradient Descent Algorithm}

Adaptive algorithms such as AdaGrad \citep{Duchi_AdaGrad_2011}, RMSProp \citep{Tieleman_RMSP_2012}, and Adam \citep{Diederik_Adam_2017} have demonstrated their effectiveness. These methods aim to automatically rescale each direction of the gradient descent by using an adaptive learning rate, thereby adjusting the scale of individual gradient components and ensuring a more balanced evolution of the optimization process.

In our work, we choose to leverage the advantages of the AdaGrad algorithm. Its primary strength lies in the adaptivity of its learning rate, which is crucial since the convergence behavior of stochastic gradient algorithms is highly sensitive to this choice. Additionally, AdaGrad’s component-wise independence in adjusting the learning rate proves particularly beneficial when handling the non-differentiable penalty term involved in the criterion discussed in the next section.

The algorithm proceeds in two main steps. First, we simulate a realization of the latent variables, either by directly sampling from their posterior distribution using the current parameter values or by performing a step of a Metropolis-Hastings sampler \citep{marin2007bayesian}. In the second step, the parameter estimates are updated via adaptive  preconditioned gradient descent.  More precisely, given the sampled latent variables, we compute the gradient of the complete log-likelihood and, for each parameter component, cumulate the squared gradients over iterations. These cumulative values are then used to define a preconditioning diagonal matrix, with entries inversely proportional to the square root of the cumulated gradients.

 \subsection{Adaptive  Weighted Proximal Stochastic Gradient Descent Algorithm} 
\label{sec:proximalweighted}

     Due to the non-differentiability of the considered $\ell_1$-penalty term, we  use a proximal algorithm as presented by \cite{achab_learning_2017} and \cite{fort_stochastic_2017}. Thus, we add a proximal step in the adaptive stochastic  gradient algorithm to calculate the penalized estimate defined in Equation  \eqref{eq:penalized_estimator}.

        The algorithm is divided into three steps;  a realization of the latent variables is sampled with a first step called \textit{Simulation}, which uses a Metropolis-Hastings sampler (see \cite{marin2007bayesian}). The second step is the classical gradient descent on the gradient of the logarithm of the  complete likelihood, the \textit{Forward} step. 
    The last step, called \textit{Backward}, deals with the penalty term. We apply the classical proximal operator (see \cite{moreau_fonctions_1962,rockafellar_monotone_1976}), defined below :

        \begin{equation}\label{eq:prox}
            \prox_{\pen_\regularization}(\beta) = \argmin_{\beta' \in \setr^p}
            \left(\pen_\regularization(\beta') + \frac 12\norm{\beta-\beta'}_2^2 \right).
        \end{equation}

        Considering the $\ell_1$-penalty, the proximal operator has an explicit form given for all $ i \in   \{1, ..., p \}$:
        \begin{equation}\label{eq:prox_lasso}
            (\prox_{\regularization | \ \cdot \ |}(\beta))_i = \stackeq[cr]{
                0 &\text{if } \abs{\beta_i}<\regularization 
            \\ \beta_i - \regularization &\text{if } \beta_i \geq \regularization
            \\ \beta_i + \regularization &\text{if } \beta_i \leq -\regularization} 
        \end{equation} 
        
        The \textit{Backward} step corresponds to the application of the proximal operator on the result of the \textit{Forward} step. 
        
        \begin{algorithm}[!ht] \label{algo:SPG-fim}
                \caption{Adaptive  Weighted Proximal Stochastic Gradient Algorithm}
                \Require{$\step_0>0$, $\epsilon>0$, $K_{max} \in \mathbb{N}$ } 
                \Initialize Initialize starting point $\theta_0 \in \mathbb{R}^d$; adaptive stepsize vector $s_0 =0$; iteration $k=0$.
                
                \While{$(\theta_\iter)$ not converged or $k \leq K_{max} $}{
                    $\iter \leftarrow \iter +1$
                    \\
                    \myBullet {\bf Simulation step:} \\ 
                    \quad Draw $\lat^{(\iter)}$
                    from the posterior distribution $p( . |  \theta_k)$ or
                    using a single step of a  Metropolis Hastings algorithm  having as stationary distribution the  posterior.
                    \\
                    \myBullet {\bf Gradient computation:}
                    Compute $v_k = \frac 1\nobs \som i1\nobs \nabla \log f_{\theta_k}(\varobs[_i],\varlat[^{(k)}_i])$
                    \\
                    \myBullet {\bf Adaptive stepsize:}
                    $s_{\iter,l}^2 \leftarrow s_{\iter-1,l}^2+ v_{\iter,l}^2$ for $l \in\{1,\dots, \dimcovgd\}$
                    \\
                    \myBullet {\bf Update parameters:} \\
                    \quad \myBullet {\bf Forward step:} 
                    $\omega_{\iter,l} \leftarrow \theta_{\iter-1,l} + \step_0 v_{\iter,l}/(s_{\iter,l}+\epsilon)$ for $l \in\{1,\dots, \dimcovgd\}$
                    \\
                    \quad \myBullet {\bf Backward step:} 
                    $\theta_{\iter,l} \leftarrow \prox_{(\Prec_{\iter,l}+\epsilon)/\step_0, \regularization|.|}(\omega_{\iter,l})$     for $l \in\{1,\dots, \dimcovgd\}$
                    
                }
                \Return{$\hat \theta = \theta_\iter$}
            \end{algorithm}

\section{Numerical experiments}\label{sec:simu}
    In this section, we  study the performance of the proposed procedure called \NLJMLASSO.
  We assess its performance in different scenarios. We  first study the impact of the  censoring rate on the variable selection and estimation of the parameters when individual parameters of the mixed model are independent. Then, we will highlight the difficulty appearing  when the individual parameters of the model are correlated.  
  In each scenario, we generate $n_{\text{run}}$  independent data sets and for each run  we  fit the considered model using the routine described in Section \ref{algo:metho} and get an estimate $\hat\theta^{\text{MLE}}_{i}$ for each $i \in \{1,..., n_{\text{run}}\}$.   
   To compare the results of the different scenarios, we
    evaluate the variable selection capacity of the method by computing the sensitivity, specificity  and accuracy. Sensitivity is the proportion of true positives (TP) correctly identified. Specificity is the proportion of the true negatives (TN) correctly identified. Finally, accuracy is the proportion of true results, either true positive or true negative. We abbreviated false negatives and false positives  respectively by FN and FP.

    \begin{equation*}
        \text{Se} = \frac{TP}{TP + FN} ~~
        \text{Sp} = \frac{TN}{TN + FP} ~~
        \text{Ac} = \frac{TP + TN}{P+N}
    \end{equation*}

 We   look also at the mean square error (mse) and  the relative root mean square errors (RRMSE) to measure the estimation quality of the method defined as:
    \begin{equation*}
        \text{RRMSE}(\hat\theta^{\text{MLE}}_{k}) = \sqrt{\frac 1{n_{\text{run}}}\sum_{i=1}^{n_{\text{run}}} \frac{(\hat\theta^{\text{MLE}}_{i,k} - \theta_k^*)^2}{\theta_k^{*2}}} ~~
        \text{mse}(\hat\theta^{\text{MLE}}_{k}) = \frac 1{n_{\text{run}}}\sum_{i=1}^{n_{\text{run}}} (\hat\theta^{\text{MLE}}_{i,k} - \theta_k^*)^2
    \end{equation*}

where $\theta^*$ denotes the true value of the parameters.

        \subsection{Variable selection in high dimension and re-estimation analysis in a non linear joint model with  independent random effects from censored survival data} 
        \renewcommand{\varlat}{\omega} 
        \renewcommand{\vareps}{\sigma}

        We consider the joint causal model defined in Section \ref{chap4:sec:model} with the classical logistic growth function for the non linear mixed-effect model defined by:
        \begin{equation}\label{chap4:eq:logistique}
            \obs_{i,j} 
                = \frac{\lat_{i1}}{1+\exp\left(-\left(\covrep_{i,j} - \lat_{i2}\right) / \ldp\right)} +\eps_{i,j}
        \end{equation}
        \noindent where the quantity $\lat_{i1}$ represents the asymptotic  value of the curve, $\lat_{i2}$  the value of the sigmoid midpoint, and $\ldp$  the logistic growth rate. The term $\eps_{i,j}$ represents an  additive noise term. We model for each individual $i$ the corresponding individual parameter $\lat_i\in\setr^2$ through a Gaussian random variable with expectation $\meanlat=(\meanlat_1,\meanlat_2) \in \setr^2$ and diagonal covariance matrix $\Omega = diag(\varlat_1^2,\varlat_2^2)$ with $\varlat_1^2,\varlat_2^2 \in(\setr_+^*)$. We assume that the noise term follows a centered Gaussian distribution 
with unknown variance $\vareps^2 \in\setr^*_+$.
        Regarding the survival model, we consider a fixed unknown constant value for the baseline $h_0(t) = \hazard_0, \forall t\geq 0$. Moreover we assume that the time of interest $\mbf T^*$ can be related to the  longitudinal dynamic using the following link function: $f(\alpha,\mc M(t;\lat_i)) = \alpha m(t,\lat_i)$ with
         \begin{equation}\label{chap4:eq:mfunction}
       m(t,\lat_i)
                = \frac{\lat_{i1}}{1+\exp\left(-\left(\covrep_{i,j} - \lat_{i2}\right) / \ldp\right)} , \  \forall t \geq 0
        \end{equation}

      Thus the  hazard rate writes: $\forall i \in\{1,...,\nobs\}, \forall j\in \{1,...,\nrep\}$ :
        \begin{equation}\label{chap4:eq:joint_simu}
            \hazard(t|\mc M(t,\lat_i), \covgd_i) 
                = \hazard_0 \exp \bigr(\beta^T \covgd_i + \alpha \memfct(t, \lat_i) \bigl), \  \forall t \geq 0
        \end{equation}
        The whole joint model parameter are denoted by $\theta = (\hazard_0, \beta, \alpha, \meanlat_1, \meanlat_2, \ldp, \varlat_1^2, \varlat_2^2, \vareps^2)$. We generated independently $100$ data sets with the following parameter's values: $\hazard_0 = 0.05$, $\meanlat_1 = 4$, $\meanlat_2 = 15$, $\ldp = 4$, $\varlat_1^2 = 0.1$, $\varlat_2^2 = 1$, $\vareps^2 = 0.05$, and $\alpha = 1.5$.  Each  vector of covariates $\covgd_i$ is generated with $p$ components, where each components is drawn independently from a binomial distribution with a success probability of 0.2, $\covgd_{i,k}\sim \mc B(0.2), \forall i\in\{1\dots,\nobs\}, k\in\{1\dots,\dimcovgd\}$. The first four components of the vector $\beta$ are set to $(-2,-1,2,2)$, while the remaining components are set to zero. We simulate censoring of the event $\event^*_i$ by generating a random variable $C_i$ from a uniform distribution on the interval $[0, \event^*_i]$ and a censoring indicator $\delta_i$ according to an independent Bernoulli distribution with a probability $\rho$ of being observed equal to the desired observation percentage. Each individual's survival time has then a probability of being censored equal to $1-\rho$. The event time $\event_i$ is defined by $\event_i = C_i \un_{\delta_i = 0} + \event^*_i \un_{\delta_i = 1}$.  The number of individuals is fixed at $\nobs = 200$ and the number of observations per individual is fixed at $\nrep = 15$. The number of covariates is fixed at $p = 500$.        In this simulation study, we focus on the variable selection  procedure as well as on the      parameter estimation in the reduced model  once the covariate support has been selected    in a context of censored data.  Thus we consider three scenarios with different percentages of censored data : without censoring $0\%$,  with moderate censoring $40\%$, and heavy censoring $70\%$.

        The results of the selection procedure are presented  in Table \ref{chap4:tab:SeSpAc_logistic}.  The selection score are given,  as well as  the \textit{mean square error} (mse).  We observe that the censoring has a strong impact on the selection procedure. The sensitivity, which is the proportion of relevant covariates selected, decreases strongly as the percentage of censoring increases.  Figure  \ref{chap4:fig:selection_percentage_censoring_logistic}
        presents the proportion of correctly selected non-zero components for the first four components of $\beta$. We observe that these proportions decrease strongly as the censoring rate increases. In all scenario the procedure struggles to select non-zero components of $\beta$ when the censoring is high. The component $\beta_2$  seems to be more difficult to select than the others. This can be explained by its relative lower amplitude.

        \begin{tableht}
            \caption{Average sensitivity (Se), specificity (Sp), accuracy (Ac), and mean square error (mse) using the procedure \NLJMLASSO over $100$ repetitions for the logistic joint model.}
            \label{chap4:tab:SeSpAc_logistic}
            \begin{tabular}{lccc}
                \toprule
                    & 0\% & 40\% & 70\% \\
                \midrule
                Ac & 0.999 & 0.997 & 0.994 \\ 
                Se & 0.940 & 0.782 & 0.230 \\ 
                Sp & 1.000 & 0.999 & 1.000 \\ 
                mse & 0.001 & 0.005 & 0.019 \\ 
                \bottomrule
                \end{tabular}
        \end{tableht}\vspace*{3mm}

        \begin{figureht}
            \includegraphics[width=0.7\textwidth]{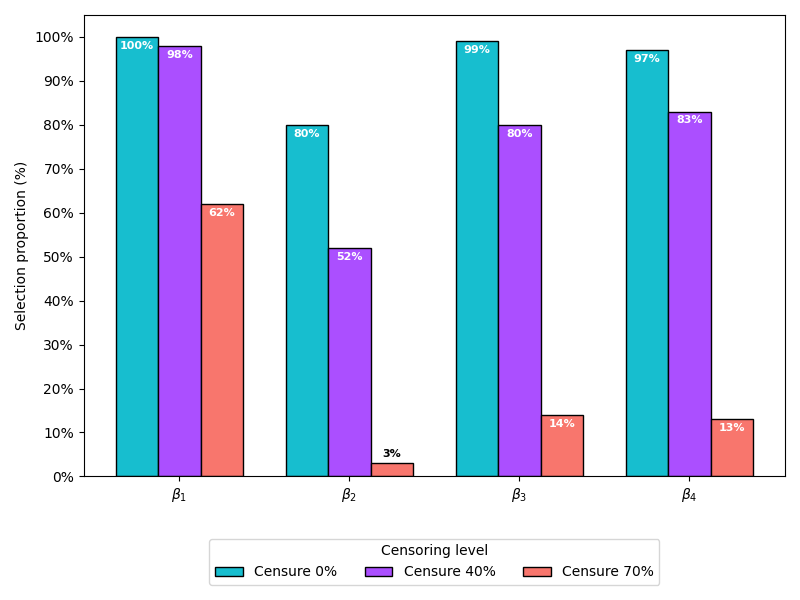}
            \caption{Proportion of correctly selected non-zero components for the logistic  joint model using the \NLJMLASSO procedure over $100$ repetitions with three different censoring rates.} \label{chap4:fig:selection_percentage_censoring_logistic}            
        \end{figureht}

        After the selection procedure, we re-estimate the parameters in the reduced model. The results are summarized in Table  \ref{chap4:tab:RMSE_logistic}. The censoring rate seems to have globally a small  impact on the estimation of the parameter of the mixed-effect models in this configuration.  The variances of the random effects are poorly estimated whatever the censoring rate.  The variance $\varlat_1^2$ seems to be more affected than the other parameters. The link parameter $\alpha$ is more difficult to estimate when the censoring  increases (from $10\%$ to $16\%$ of relative error) as well as the baseline constant $\hazard_0$.  The four non-zero components are well estimated without censoring. However these estimations are deeply deteriorated under moderate and heavy censoring, in particular for the component $\beta_2$. This phenomenon is coherent with results presented in Figure \ref{chap4:fig:selection_percentage_censoring_logistic} and can be explained by the lower proportion of correctly selected variables for this specific component.

        \begin{tableht}
            \caption{Mean  parameter estimates and average relative root mean square errors (RRMSE)  over $100$ repetitions for the logistic joint model  where $\hat\theta$ is the estimates in the reduced model after variable selection using the adaptive stochastic gradient algorithm.  The value $\para^*$ denotes  the true value of the parameter.}
            \label{chap4:tab:RMSE_logistic}
            \resizebox{0.7\textwidth}{!}{
            \begin{tabular}{@{}cc*{4}{cc}@{}} 
                \toprule
                 & & & \multicolumn{2}{c}{ C = 0\% }  & \multicolumn{2}{c}{ C = 40\% }  & \multicolumn{2}{c}{ C = 70\% } \\ \cmidrule(r){4-5} \cmidrule(r){6-7} \cmidrule(r){8-9}
                & & $\para^*$ & $\hat\theta$ & RRMSE($\%$) & $\hat\theta$ & RRMSE($\%$) & $\hat\theta$ & RRMSE($\%$)  \\
                \midrule
    & $\meanlat_1$ & 4.0 & 4.15 & 8.9 & 4.14 & 10.2 & 4.14 & 15.0 \\ 
    & $\meanlat_2$ & 15.0 & 15.25 & 4.5 & 15.23 & 5.3 & 15.21 & 7.2 \\ 
    & $\varlat_1^2$ & 0.1 & 0.226 & 152.2 & 0.255 & 177.9 & 0.282 & 204.8 \\ 
    & $\varlat_2^2$ & 1.0 & 0.69 & 39.6 & 0.67 & 44.0 & 0.67 & 44.5 \\ 
    & $\ldp$ & 4.0 & 4.06 & 4.3 & 4.05 & 4.7 & 4.04 & 5.5 \\ 
    & $\vareps^2$ & 0.0 & 0.050 & 6.8 & 0.050 & 7.3 & 0.049 & 8.1 \\ 
    & $\alpha$ & 1.5 & 1.49 & 14.3 & 1.59 & 20.0 & 1.29 & 26.9 \\ 
    & $h_0$ & 0.05 & 0.051 & 17.5 & 0.036 & 32.0 & 0.024 & 53.3 \\ 
    & $\beta_{1}$ & -2.0 & -2.01 & 15.5 & -1.98 & 25.2 & -1.19 & 64.2 \\ 
    & $\beta_{2}$ & -1.0 & -0.87 & 48.2 & -0.63 & 72.9 & -0.04 & 99.0 \\ 
    & $\beta_{3}$ & 2.0 & 2.00 & 17.3 & 1.64 & 46.7 & 0.33 & 93.1 \\ 
    & $\beta_{4}$ & 2.0 & 1.95 & 19.7 & 1.69 & 43.3 & 0.31 & 93.7 \\ 
                \bottomrule
            \end{tabular}}
        \end{tableht} \vspace*{3mm}

    As an illustration of the convergence of the adaptive stochastic gradient  algorithm,   Figure \ref{chap4:fig:theta_iter_example} displays five  trajectories  of   the re-estimation of  survival parameters $(\hazard_0, \alpha, \beta)$ as a function of iterations during the execution with different initializations.
     We observe that the algorithm converges to a limit close to the true value of the parameters. We emphasize here that the value targeted by the adaptive stochastic gradient algorithm is a critical point of the marginal likelihood,  which  turns out to be a maximum of the marginal likelihood.  We recall that the asymptotic behavior expected under regularity assumptions as the sample size increases is indeed a convergence of the maximum likelihood estimate towards the true value of the parameter. This explains the fact that the limits are not exactly equal to the true values  of the parameters,  but only close to them.

        \begin{figureht}
            \includegraphics[width=0.7\textwidth]{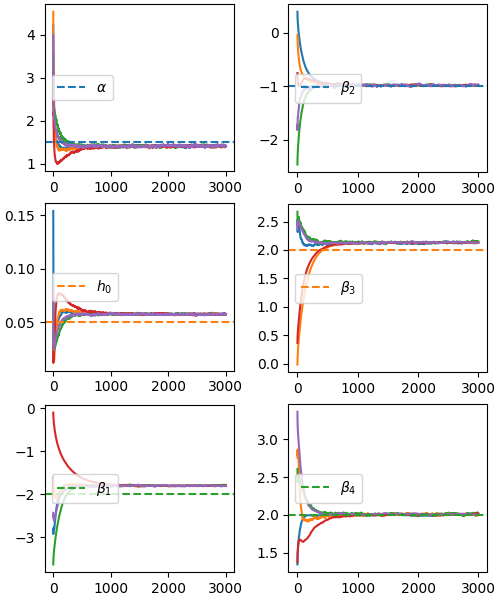}
            \caption{Five trajectories of the adaptive stochastic gradient algorithm  with different initializations for the re-estimated survival parameters in the logistic joint model without censoring. The true value of the parameters is represented by a dashed line.} \label{chap4:fig:theta_iter_example}
        \end{figureht}

        \subsection{Variable selection in high dimension and re-estimation analysis in a non linear joint model with   censored survival data in the presence of correlated individuals parameter}

        In this section, we investigate the effect of the presence of a correlation between the individual parameters on the variable selection and re-estimation for the joint model defined in the previous section.   For this purpose, we consider that  the bivariate individual parameters $\lat_i$  are drawn from a correlated Gaussian random variables with a covariance matrix $\Omega = \begin{pmatrix} \varlat^2_1 & \varlat_{12} \\
        \varlat_{12} & \varlat^2_2 \end{pmatrix} \in\setsppr{2 }$. We set $\varlat^2_1 = 0.1$, $\varlat^2_2 = 1$, and $\varlat_{12} = 0.1$, which corresponds to  a correlation coefficient of $0.31$ between the two individual parameters. The other parameters are set as in the previous section.
     We simulate independently $100$ data sets with the same censoring procedure as previously and with the same censoring rates $0\%$, $40\%$ and $70\%$. The results of the variable selection are given in Table  \ref{chap4:tab:SeSpAc_logisticb}.  We observe the same  impact of the censoring percentage on the quality of the selection.  However Figure  \ref{chap4:fig:selection_percentage_censoring_logistic_b} shows that the different non zero components of $\beta$ are selected less accuratly, since we observe in mean a loss of   $10\%$ in the  proportion of correct selected variables  compared to the previous model without correlation. 

        \begin{tableht}
            \caption{Average sensitivity (Se), specificity (Sp), accuracy (Ac), and mean square error (mse) using the procedure \NLJMLASSO over $100$ repetitions for the logistic joint model with correlated individuals’ parameters.}
            \label{chap4:tab:SeSpAc_logisticb}
            \begin{tabular}{lccc}
                \toprule
                    & 0\% & 40\% & 70\% \\
                \midrule
                Ac & 0.999 & 0.997 & 0.994 \\ 
                Se & 0.944 & 0.741 & 0.247 \\ 
                Sp & 1.000 & 0.999 & 1.000 \\ 
                mse & 0.001 & 0.006 & 0.019 \\ 
                \bottomrule
                \end{tabular}
        \end{tableht}\vspace*{3mm}

        \begin{figureht}
            \includegraphics[width=0.7\textwidth]{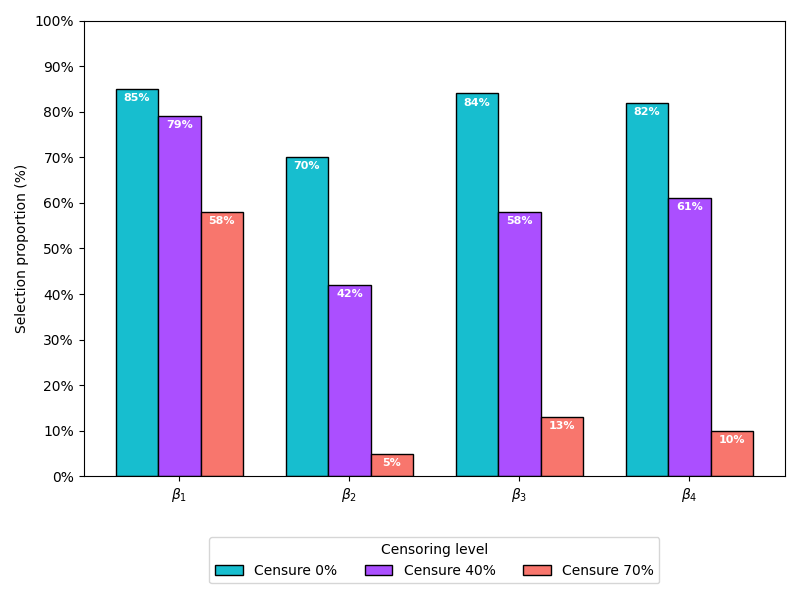}
            \caption{Proportion of correctly selected non-zero components for the logistic  joint model with correlated individuals’ parameters  using the \NLJMLASSO procedure over $100$ repetitions for three different censoring rates.} \label{chap4:fig:selection_percentage_censoring_logistic_b}            
        \end{figureht}

        The results of the re-estimation of the parameters after the selection procedure are presented  in Table \ref{chap4:tab:RMSE_logisticb}. The presence of correlation between the individual parameters seems to have  a strong impact on the estimation of the covariance matrix $\Omega$.  Estimations of other 
         parameters seem to be  less  affected by the presence of a correlation and have similar estimation results as in the previous section. 

        \begin{tableht}
            \caption{Mean parameter estimates  and average relative root mean square errors (RRMSE) over $100$ repetitions for the logistic joint model  with correlated individual's parameters, where $\hat\theta$ are the estimates in the reduced model after variable selection using the adaptive stochastic gradient algorithm. The value  $\para^*$ is the true value of the parameter.}
            \label{chap4:tab:RMSE_logisticb}

        \resizebox{0.8\textwidth}{!}{
            \begin{tabular}{@{}cc*{4}{cc}@{}} 
                \toprule
                 & & & \multicolumn{2}{c}{ C = 0\% }  & \multicolumn{2}{c}{ C = 40\% }  & \multicolumn{2}{c}{ C = 70\% } \\ \cline{3-4} \cline{5-6} \cline{7-8}
                & & $\para^*$ & $\hat\theta$ & RRMSE($\%$) & $\hat\theta$ & RRMSE($\%$) & $\hat\theta$ & RRMSE($\%$)  \\
                \midrule
    & $\meanlat_1$ & 4.0 & 4.19 & 10.6 & 4.26 & 14.1 & 4.21 & 16.5 \\ 
   & $\meanlat_2$ & 15.0 & 15.261 & 4.9 & 15.380 & 6.5 & 15.26 & 8.1 \\ 
   & $\varlat_1^2$ & 0.1 & 0.36 & 342.5 & 0.42 & 453.1 & 0.55 & 555.1 \\ 
   & $\varlat_{12}$ & 0.1 & 0.372 & 388.4 & 0.39 & 459.8 & 0.56 & 599.5 \\ 
   & $\varlat_2^2$ & 1.0 & 1.27 & 67.0 & 1.34 & 90.5 & 1.67 & 120.7 \\ 
   & $\ldp$ & 4.0 & 4.05 & 4.0 & 4.07 & 5.1 & 4.04 & 6.4 \\ 
   & $\vareps^2$ & 0.0 & 0.05 & 6.8 & 0.05 & 7.6 & 0.049 & 7.7 \\ 
   & $\alpha$ & 1.5 & 1.50 & 14.8 & 1.53 & 20.0 & 1.31 & 25.8 \\ 
   & $h_0$ & 0.05 & 0.05 & 17.7 & 0.037 & 30.2 & 0.02 & 53.8 \\ 
   & $\beta_{1}$ & -2.0 & -2.01 & 16.4 & -1.90 & 25.1 & -1.235 & 60.9 \\ 
   & $\beta_{2}$ & -1.0 & -0.895 & 45.5 & -0.63 & 72.7 & -0.09 & 98.4 \\ 
   & $\beta_{3}$ & 2.0 & 2.00 & 17.1 & 1.478 & 55.1 & 0.346 & 92.7 \\ 
   & $\beta_{4}$ & 2.0 & 1.94 & 20.9 & 1.54 & 51.1 & 0.29 & 94.5 \\ 
                \bottomrule
            \end{tabular}}
        \end{tableht} \vspace*{3mm}

\section{Conclusion and perspectives}

    In this work, we jointly addressed variable selection and parameter estimation in a  non linear joint model for survival data and longitudinal outcomes. We connected a survival model with a non linear mixed-effects model through a link function and covariates  in high dimension.  To perform the variable selection procedure in the joint model, we consider as criterion the marginal log-likelihood reguralized with a $\ell_1$-penalty term and 
   run an adaptive  weighted proximal stochastic gradient algorithm  to deal simultaneously with the latent variables and the non-differentiability of the  penalty.   Then we rely on an eBIC criterion to choose an optimal value for the regularization parameter. We   finally re-estimate the parameters in the reduced model with an adaptive stochastic gradient descent algorithm.   Our methodology has been thoroughly evaluated on simulated data to demonstrate its performance. We evaluate the impact of censoring and of the presence of correlation between individual parameters.

    One interesting perspective of this work consists of addressing the prediction task carefully. It would be interesting to set up a method for predicting survival time from some start of longitudinal data observation. Previous studies, such as \cite{desmee_nonlinear_2017} and \cite{kerioui_modelling_2022}, have explored this task within the context of joint models without high dimensional covariates.
   Furthermore, considering the partial likelihood for the survival part rather than the parametric complete  likelihood would be an interesting  avenue to explore. This approach would eliminate the need to make assumptions about the form of the baseline risk, which is especially relevant when dealing with the Cox model. Indeed, when considering the Cox model, the partial likelihood allows for estimating consistently the regression parameter without knowing the baseline function. It would be interesting to adapt the proposed procedure to the 
    partial likelihood instead of the parametric complete  one. 

\section*{Funding and Acknowledgements}

     This work was funded by the https://stat4plant.mathnum.inrae.fr/(Stat4Plant) project ANR-20-CE45-0012.


    \printbibliography

\section*{Appendix}

    \begin{algorithm}
        \caption{Simulation Algorithm}
        \Require{number of individuals N, number of longitudinal observation $J$, parameter values $\sigma^2, \beta$} 

        \myBullet \bf{Generate :} \\ 
            \begin{itemize}
                \item Set $t_j = t_max + j \times \frac{t_{max}-t_{min}}J ~\forall j \in \{1, ..., J\}$
                
                \item \bf{Generate longitudinal observation :}
                \begin{itemize}
                    \item draw $\epsilon \under\sim{i.i.d.} \mc N(0, \sigma^2)$,
                    \item draw $\lat_i \under\sim{i.i.d.} \mc N(\mu, \Gamma) $,
                    \item compute $\mbf Y_{i,j} = \memfct(t_j; \lat_i) + \varepsilon_{i,j}$.
                \end{itemize}

                \item  \bf{Generate covariates :}
                \begin{itemize}
                    \item draw $U_i^* \sim \mc U([-1,1])$,
                    \item define $U_i = U_i^* - \bar U_i^*$, where $\bar U_i^*$ is the mean vector along the columns.
                \end{itemize}

                \item  \bf{Generate Survival times : $\forall i \in \{1, ..., N\}$}
                \begin{itemize}
                    \item generate covariates $U_i^* \sim \mc U([-1,1])$,
                    \item draw $u_i \sim U([0,-1])$,
                    \item solve : $T_i^* = F_i^{-1}(u_i)$ where $F_i(t) = 1-S_i(t)$ where $S_i$ is the survival function $S_i(t) = \exp\left(-\int_0^t h_i(t|\mc M(t,\lat_i)dt \right)$
                \end{itemize}

                \item \bf{Censoring data}
                \begin{itemize}
                    \item set $C_{i}$ such that $\frac{\#\{ T_i| T_i \leq C_{i} \}}N$ is the right percentage of censoring
                    \item censor by setting $T_i = \min(T_i^*, C_{i})$
                    \item compute $\delta_i = \un_{T_i^* \leq C_{i}}$
                    \item remove all longitudinal observations after the time $C_{i}$
                \end{itemize}
                
            \end{itemize}
        
    \end{algorithm}


\label{lastpage}

\end{document}